# THE FLOYD-WARSHALL ALGORITHM, THE AP AND THE TSP

**Howard Kleiman**

## 1. INTRODUCTION.

**Let s and t be two vertices of a connected weighted graph G represented by the matrix $M_G$. The *shortest path problem* finds a path between s and t whose total edge weight is minimum. Generally, edge-weight is taken to mean *distance* but the word is used loosely and may represent some other measurable quantity. Dijkstra's algorithm [1] finds the distance between s and *all* of the other vertices of G. However, it assumes that all edge weights are non-negative. The Floyd-Warshall algorithm [2], [7], finds the shortest paths between *all pairs* of nodes. Furthermore, unlike the Dijkstra algorithm, it allows arc weights to be negative. In this paper, we use the variant of the F-W algorithm given in [6]. Of particular interest to us, this algorithm allows us to find all cycles of smallest total arc-weight. The running time of the F-W algorithm is $O(n^3)$. In this paper, given a random n-cycle, D, obtained from the symmetric group, $S_n$, we apply $D^{-1}$ to permute the columns of $M_G$. Thus, if the arc (a, D(b)) has the weight w(a, D(b)) in $M_G$, its weight is transformed into w(a, b) in $D^{-1}M_G$. It follows that a negative- weighted cycle, C, in $D^{-1}M_G$, say $(a_1 a_2 \ldots a_i)$, has the property that the total edge weight of DC is less than that of D. Furthermore, as long as no pair of consecutive vertices of C - $a_j a_{j+1}$ - has the property that $a_{j-1} a_j$ is an arc of D, DC will be a *derangement,* i.e., a permutation that moves all points in V = {1,2,...,n}. In [4], we used *H-admissible* permutations to obtain hamilton circuits in graphs. In the algorithm given here, we use *admissible permutations* as defined in [5]. Starting with a random n-cycle, they were used in [5] to obtain a spanning set of |log n| disjoint cycles. Here they are used to obtain a rough approximation to an optimal solution of the edge-distance assignment problem defined by the entries in $M_G$. This approximation yields a weight matrix in which the negative entries are generally fewer in number and smaller in absolute value than those in $M_G$. Theorems 1 and 2 allow us to use considerably fewer trials than generally used in Floyd-Warshall. Once we have obtained an optimal solution, $\sigma_{APOPT}$, to the assignment problem, the only cycles we can obtain in $\sigma_{APOPT}^{-1} M_G$ are all positive. We then use the above theorems and the F-W algorithm to obtain the smallest positive cycles of total edge weight less than a fixed positive number, say N. As we proceed, we test cycles and sets of disjoint cycles to see if any product is an n-cycle. The n-cycle of smallest weight that we obtain is an approximate optimal solution** **to the TSP. The reason that it is not exact is that the F-W algorithm allows us to obtain cycles of *smallest* weight constructed by always using the shortest distance between any pair of vertices. It is possible that a cycle, C, may be a disjoint cycle of a permutation, s, such that $\sigma_{APOPT} s = \sigma_{TSPOPT}$ where C contains a subpath [a, ..., b] which is not the shortest path between a and b. However, we do give conditions where our version of the F-W algorithm leads to $\sigma_{TSP}$. Furthermore, we give necessary conditions in Theorem 6 that cycles of s not obtainable by F-W must satisfy.**

## II. THEOREMS

**THEOREM 1. Let C be a cycle of length n. Assume that the weight w (i, C(i)) (i=1,2,...,n) corresponds to the arc (i, C(i)) of C. Then if**

$$W = \sum_{i=1}^{i=n} w(i, C(i)) < 0$$

**there exists at least one value i = i' with $1 \le i_1 \le n$ such that**

$$\sum_{j=0}^{j=m} w(i'+j, C(i'+j)) < 0 \qquad \text{(A)}$$

**where m = 0,1,2,...,n-1 and $i' + j > n$ is represented by its value modulo n.**

**PROOF. The proof is by induction. Let n = 2. Then one of two cases occurs. Either both arcs have negative weights or one is negative and one is positive. In the latter case, since the sum is negative, if our first arc is negative, the theorem is proved. Assume that the theorem is true when the cycle has n arcs. We now consider the case when it has n+1 arcs. In this case, there must exist at least one case, say $i = i^*$, where $w(i^*) + w(i^* + 1) < 0$. Suppose this wasn't the case. Then as we traversed the cycle in a clockwise manner, the sum of the weights of any two consecutive arcs would always be non-negative. Therefore, the total sum of the weights of the cycle couldn't be negative. We now replace the nodes $i^*$ and $i^* + 1$ by the node $i^{**}$ whose weight is defined as $w(i^{**}) = w(i^*) + w(i^* + 1)$. The cycle C having n+1 nodes has now been replaced by the cycle C* having n nodes. By induction, the theorem is valid for C*. Let $i'$ be the node of C* defined in Theorem 1 and $i' + j^* = i^{**}$. Then**

$$\sum_{j=0}^{j=j^*} w(i' + j, C(i' + j)) < 0$$

**Now consider the sum**

$$\sum_{j=0}^{j=j^*-1} w(i' + j, C^*(i' + j)) + w(i^*, C(i^*))$$

**The summation on the left is negative due to our inductive assumption that the theorem holds for C*. On the other hand, we have assumed that** $W = \sum_{i=1}^{i=n} w(i, C(i)) < N$ and $w(i^*, C(i^*))$ is negative. Therefore, the sum of these two quantities is also negative. Furthermore,

$$\sum_{j=0}^{j=j^*} w(i'+j, C(i'+j)) = \sum_{j=0}^{j=j^*-1} w(i'+j, C(i'+j)) + w(i^*, C(i^*)) + w(i^*+1, C(i^*+1)) < 0$$

Thus, if we now go back to C, using the same node $i'$, we can replace the node $i^{**}$ of C* by the nodes $i^*$ and $i^*+1$ of C and still satisfy condition A of the theorem, concluding the proof by induction.

**COROLLARY.** Suppose that C is a cycle such that

$$W = \sum_{i=1}^{i=n} w(i, C(i)) < N \qquad \text{(B)}$$

Then there exists a value $i = i'$ such that each partial sum, $S_j$, has the property that

$$S_j = \sum_{j=0}^{j=m} w(i'+j, C(i'+j)) < N \qquad \text{©}$$

always holds.

**PROOF.** Subtract N from both sides of (B). Now let the weight of arc (n, C(n)) become w(n, C(n)) – N. From Theorem 1,

$$W^* = W - N = \sum_{i=1}^{i=n} w(i, C(i)) < 0 \quad \text{(D)}$$

has the property that every partial sum is less than 0 where each partial sum is of the form given in C. It follows that if we add N to each side of D and restore w(n, C(n)) to its original value, every partial sum is less than N.

**EXAMPLE.** We now give an example of how to obtain $i = i'$.

Let C = ($a_1$ $a_2$ ... $a_n$),

$a_1$ = -7, $a_2$ = -10, $a_3$ = +1, $a_4$ = +2, $a_5$ = -7, $a_6$ = +4, $a_7$ = -9, $a_8$ = +11, $a_9$ = -2, $a_{10}$ = -1, $a_{11}$ = -4, $a_{12}$ = -4, $a_{13}$ = -8, $a_{14}$ = +9, $a_{15}$ = +9, $a_{16}$ = +21, $a_{17}$ = +1, $a_{18}$ = -2, $a_{19}$ = -1, $a_{20}$ = -3, $a_{21}$ = -3, $a_{22}$ = -12, $a_{23}$ = +6, $a_{24}$ = +2, $a_{25}$ = +3

-7 –10 +1 +2 –7 +4 –9 +11 –2 –1 –4 –4 –8 +9 +9 + 21 +1
-2 -1 -3 -3 -12 +6 +2 +3

We now add terms with like signs going from left to right. We place the ordinal number of the first number in each sum above it.

| | | | | | | | | | |
|---|---|---|---|---|---|---|---|---|---|
| *1* | *3* | *5* | *6* | *7* | *8* | *9* | *14* | *18* | *23* |
| -17 | +3 | -7 | +4 | -9 | +11 | -19 | +40 | -21 | +11 |

We next add the positive number to the right of each negative number to the negative number.

| | | | | |
|---|---|---|---|---|
| *1* | *5* | | | *18* |
| -14 | -3 | +2 | +21 | -10 |

We now add terms with like signs going from left to right.

| | | |
|---|---|---|
| *1* | | *18* |
| -17 | +23 | -10 |

**Finally, assuming that all points lie on a circle, we add like terms going from left to right. We thus obtain**

***18*** ***18***

**-27 +23 = -4**

**This tells us that i' is the eighteenth ordinal number, namely, –2.**
**Thus, the partial sums are: -2, -3, -6, -9, -21, -15, -13, -10, -17, -27, -26, -24, -31, -27, -36, -25, -27, -28, -32 –36, -44, -35, -26, -5, -4**
**The following theorem is trivial to prove but essential for the use of MIN(M) throughout the algorithm.**

**THEOREM 2. THE FLOYD-WARSHALL ALGORITHM.**
**If we perform a triangle operation for successive values j = 1,2,...,n, each entry $d_{ik}$ of an n X n cost matrix M becomes equal to the value of the shortest path from i to k provided that M contains no negative cycles.**
**The version given here is modeled on Theorem 6.4 in [6].**
**PROOF: We shall show by induction that that after the triangle operation for $j = j_0$ is executed, $d_{ik}$ is the value of the shortest path with intermediate vertices $v \le j_0$, for all i and k. The theorem holds for $j_0 = 1$ since v = 0. Assume that the inductive hypothesis is true for $j = j_0 - 1$ and consider the triangle operation for $j = j_0$:**

$$d_{ik} := \min \{d_{ik}, d_{ij_0} + d_{j_0k}\}.$$

**If the shortest path from i through k with $v \le j_0$ doesn't pass through $j_0$, $d_{ik}$ will be unchanged by this operation, the first argument in the min-operation will be selected, and $d_{ik}$ will still satisfy the inductive hypothesis. On the other hand if the shortest path from i to k with intermediate vertices $v \le j_0$ does pass through $j_0$, $d_{ik}$ will be replaced by $d_{ij_0} + d_{j_0k}$. By the inductive hypothesis, $d_{ij_0}$ and $d_{j_0k}$ are both optimal values with intermediate vertices $v \le j_0 - 1$. Therefore, $d_{ij_0} + d_{j_0k}$ is optimal with intermediate vertices $v \le j_0$.**
**We now give an example of how F-W works.**
**Example. Let d(1, 3) = 5, d(3, 7) = -2, d(1, 7) = 25. Then d(1, 3) + d(3, 7) < d(1, 7).**
**Note however, that the intermediate vertex 3 comes from the fact that we have reached column j = 3 in the algorithm. We now substitute d(1, 3) + d(3, 7) = 3 for the entry in (1, 7). Suppose now that d(1, 10) = 7 while d(7, 10) = -5.**

$$d(1, 3) + d(3, 7) + d(7,10) = -2 < d(1, 10) = 7.$$

**The intermediate vertices are now 3 and 7. Thus, -2 can be substituted for 7 in entry (1, 10). Again, we must have reached column 7 before this substitution could be made. Now suppose that d(10, 1) = 1. Then our negative path implies that the negative cycle (1 3 7 10) exists.**

**In general, when using our specialized version of the F-W algorithm, we choose the *first* negative cycle we obtain. The reason for this becomes evident when we consider the following theorem:**

**THEOREM 3. Let M be a value matrix containing both positive and negative values. Suppose that M contains one or more negative cycles. Then if a negative path P becomes a non-simple path containing a negative cycle, N, as a subpath, C, C is obtainable as an independent cycle in the *F-W negatively-valued subpaths* algorithm (F-W n-vs algorithm as used in Phase 2).**

**Before going on the proof of the theorem, we note that subpath when we reach column $j = j_0$, N occurs as a simple cycle with $j \leq j_0$. Towards the end of Phase 2, it may require going through our set of columns from 1 to n more than once.**

**PROOF. The reason that the F-W algorithm always gives a simple path if M contains no negative cycles is because if a non-negative cycle becomes a subpath of a here-to-fore simple path, P, since the cycle doesn't *decrease* the value of the P, it is not included in either the value of P nor its arcs. On the other hand, a negative cycle *decreases* the value of P. Furthermore, each time we reach the vertex which first made the path non-simple, we would *again* decrease the value of the path. That is the reason that the version of F-W used is the one given in [6]: it doesn't yield an infinite loop if you obtain a negative cycle as a subpath. Now the proof. Suppose that we have the path**

$$P = [a_1\ a_2\ \ldots\ a_i\ v\ a_{i+2}\ \ldots\ a_{i+p}\ j_0\ \ a_{i+p+1} \ldots a_{i+p+2} \ldots a_{i+p+q}\ \ v]$$

**where $[v \ldots j_0 \ldots v]$ defines a negative cycle. Assume that we are currently constructing paths of the form $[a\ j_0]\ [\ j_0\ \ b]$ and that $v = a_{i+1}$. From our algorithm, each $a_I$, $(I = 1,2,\ldots,i+p+q)$ is less than $j_0$. Thus, the arcs in the path $[v\ a_{i+2} \ldots a_{i+p}]$ have the property that each vertex is smaller than $j_0$. The same is true for all vertices in $[j_0\ \ldots\ v]$ other than $j_0$. Since the cycle is negative, either the path $[v \ldots j_o]$ or the path $[j_0\ \ldots\ v]$ has a negative value (i.e., the sum of the values of the arcs is negative.) In the first case, v occurs before $j_0$. Thus, the negative cycle beginning with v can be found before we are finished with all possibilities for negative paths when $j = j_0$. It could occur that $v > a_1$ (I = 1,2,...i]. However, even if that is the case, the negative cycle defined by $[v \ldots j_0 \ldots v]$ will be found before we go on to $j = j_0\ + 1$.**

**The algorithm used in PHASE 2 is based on the Floyd-Warshall theorem. However, its purpose is to obtain a *negative cycle* using Theorem 1. Namely, we start with a negative-valued arc and continue with negative-valued paths until we've obtained a negative cycle.**

**THEOREM 5. Let P be a negatively-valued path in $D_h^{-1}M^{-}(k)$, say $P = [a, \ldots, i, \ldots, d, \ldots, i]$. Let i be the first vertex in P which changes P** from a simple path to a non-simple one, i.e., it is the first vertex which is repeated. Let C be the subpath of P forming the cycle (d .... i). Here d is a determining vertex of the negatively-valued cycle, C. Then we can always obtain C as an independent negative cycle in at most the next iteration, i.e., using at most n more columns in the algorithm.

PROOF. The simplest way to prove this theorem is to show that it holds in the worst possible case. Let r be the number of columns required to obtain all of the arcs in the subpath A = [a, ... , i], s, the number of columns passed through in order to obtain all of the arcs in I = [i, .... , d], t, the number of columns necessary to obtain all of the arcs in the subpath D = [d, e ... ,i]. Let P = $A \bigcup I \bigcup D$. Assume the following: $r + s + d \leq n$, all of the arcs in D other than (d, e) have non-negative values, d is a determining vertex of D, d = n, while each subsequent vertex in D is smaller than the previous one. The following example illustrates a construction for P when n = 20.

Example.

$P = [1^{-20}3^{5}7^{-5}13^{12}15^{1}19^{3}20^{-18}18^{1}14^{1}6^{3}7]$,

C = (20 18 14 6 7 13 15 19).

Note that there are no inversions in the natural order of the numbers until we reach d = 20. From that point the numbers decrease in value until they reach i = 7. From our algorithm, we can only choose an initial entry which is negative. d(1, 3) = -20, d(1, 3, 7) = -15, d(1, 3, 7, 13) = -20, d(1, 3, 7, 13, 15) = -8, d(1, 3, 7, 13, 15, 19) = -7, d(1, 3, 7, 13, 15, 19, 20) = -4. In the triangular inequalities of the algorithm, the above paths would correspond to the respective distances obtained: d(1, 3), d(1, 7), d(1, 13), d(1, 15), d(1, 19), d(1, 20). Since we obtain a sequence of negative arcs with no inversions, we can form the negative path A + I by the time we reach the column j = 20. On the other hand, starting C with the pair of arcs (20, 18), (18, 14), we obtain (20, 14) where d(20, 14) = -17, we can proceed no further until column j = 14 in $D_h^{-1}M^{-}(k+1)$ is reached where we obtain (20, 14)(14, 6) = (20, 6) with d(20, 6) = -16. The subsequent arcs in C are (6, 7), (7, 13), (13, 15), (15, 19), (19, 20). Thus, we can obtain the negative cycle C at j = 20 in $D_i^{-1}M^{-}(k+1)$. On the other hand, P reaches j = 7 in $D_h^{-1}M^{-}(k+1)$. Thus, we obtained C as a negative cycle independent of P less than n = 20 columns after P reached i. A more general proof could be presented, but I think, in this case, the example given is more transparent. The crucial idea is that P must pass through the subpaths I and D but in a different order than the way C was obtained. It is just as likely that C can be obtained independently *before* P becomes non-simple if the roles of I and D are reversed, i.e., A and I have a large number of inversions while D has none.

THEOREM 4. Let c be any real number, $S_a = \{ a_i \mid i = 1,2, \ldots ,n\}$, a set of real numbers in increasing order of value. For i = 1, 2, ..., n, let $b_i = a_i + c$. Then $S_b = \{ b_i \mid i = 1,2, \ldots ,n\}$ **preserves the ordering of $S_a$.**

**PROOF. The theorem merely states that adding a fixed number to a set ordered according to the value of its elements retains the same ordering as that of the original set.**

**Comment. This is very useful when we are dealing with entries of the value matrix which have not been changed during the algorithm. However, as the algorithm goes on, if we have entry (i, j) and d(i, j), we must go through *all* values of j = j' where the (i, j')-th entry's value in the current $\sigma_a^{-1}M^{-}(k)$ is less than the value of the (i, j')-th entry in $D^{-1}M^{-}$. (The terminology will become clearer in the examples given.)**

## III. THE ALGORITHM.

### PHASE 1.

**Let M be an n X n distance matrix. An arbitrary entry of M is d(i, j).**
**V = {1,2,...,n}**

**STEP 1. Sort each row of M in ascending order of numerical value. Call the matrix obtained MIN(M). An arbitrary entry (i, j) of MIN(M) is ORDINAL((i, j)) which stands for the ordinal value of entry (i, j) in row i.**

**STEP 2. Construct a random n-cycle, $D = (a_1\ a_2\ \ldots\ a_n)$ whose arcs are assigned the respective values of its arcs in M. Let $D^{-1}$ be the inverse of D.**

**STEP 3. Let the function ORD with respect to D be ORD($a_i$) = i. Construct $ORD^{-1}(i) = a_i$, $i = 1,2,\ldots,n$.**
**(Note that $ORD^{-1}$(i-1) = $a_{i-1}$ = $H^{-1}(a_i)$ if i $\neq$ 1. $ORD^{-1}$(1) = n.)**

**STEP 4. Given each arc ($a_i$, D($a_i$)) of D, obtain**

$$DIFF(a_i) = d(a_i, D(a_i)) - d(a_i, MIN(M)(a_i, 1)).$$

**(If MIN(M)($a_i$, 1) exists in D, we choose the second smallest value of row $a_i$, MIN(M)($a_i$, 2), etc. .**

**$MIN_j$ = min { DIFF($a_i$) $\in$ V | d($a_i$, D($a_i$)) – d($a_i$, MIN(M)($a_i$, j), j $\in$ V }**

**We then choose $MIN_2$, ..., $MIN_{|\ln n|}$, corresponding to the second, third,...,([log n] + 1)-th smallest values of DIFF($a_i$) {i =1,2,...,n}. Generally, $MIN_1$ is negative.**

**STEP 5. As we pointed out in the introduction, suppose that D = (1 7 12 17 20 15) is considered to be a permutation on the points in {1, 2 ,... , 20} . Then, given arcs (17, D(15)), (15, D(7)), (7, D(17)) of a graph G, we obtain the permutation (17 15 7) which, when applied to D, yields D* = (1 7 20 15 12 17). Since D(15) = 1, D(7) = 12, D(17) = 20, each of the given vertices appears in D*. Continuing, we first work with the value d($a_i$, MIN(M)($a_i$, 1)) obtained in $MIN_1$. We assume that it is negative. We now use Theorem 1 to obtain a cycle having a negative value such that each partial sum is negative. We note that we can *never* choose an arc of form (a, $D^{-1}$(a)): When we apply D to the terminal vertex of such an arc, we obtain (a, a), a loop. Thus, the *new permutation being constructed would not be a***

*derangement*. Continuing, we obtain $(a_i\ H^{-1}(MIN(M)(a_i, 1))$ which is given the value $d(a_i, D^{-1}(j_k))$. Letting $D^{-1}(j_k) = a_k$, we next obtain $(a_k, D^{-1}(MIN(M)(a_k, 1)) = (a_k, D^{-1}(j_m)) = (a_k, a_m)$.
If the sum of the first two terms is negative, we continue. Otherwise, we stop. If we stop, we check to see if either the sum $(a_i\ a_k) + (a_k\ a_i)$ or $(a_i\ a_k) + (a_k\ a_m) + (a_m\ a_i)$ is negative. If at least one of them is negative, we choose the one with the smaller value and save it for use later on. If the path $[a_i, a_k, a_m]$ has a negative value, we continue with this procedure. Checking the value of $(a_i\ a_k\ a_m)$, if its value is smaller than that of the first two cycles, we substitute it for the 2-cycle saved earlier. At some point in this procedure, one of the following occurs:

(1) The partial sum of the arcs of the path becomes positive.

(2) The partial sum of the arcs remains negative but the path repeats a vertex therefore forming a cycle out of a subpath of arcs.

Assume (1) occurs with the path $[a_i\ a_k\ a_m\ \ldots\ a_z]$. In this case, we save the cycle, $C_1$, with smallest negative value.
If (2) occurs, let our path be $[a_i\ a_k\ a_m \ldots a_p\ a_q\ a_r\ \ldots\ a_y\ a_m]$.
As in the first case, we have obtained the cycle, $C_2$, of smallest value from among the elements of the set

$$\{(a_i\ a_k), (a_i\ a_k\ a_m), \ldots, (a_i\ a_k\ a_m \ldots a_p), \ldots, (a_i\ a_k\ a_m \ldots a_p \ldots a_y)\}.$$

Without loss of generality, suppose the two cycles$(a_i\ a_k)$ and $(a_m \ldots a_p \ldots a_y)$ both have negative values. Call the product of the two cycles, $P_3$. Whichever of the three permutations, $C_1$, $C_2$, $P_3$ has the smallest (negative) value, gives us our best negative permutation to be applied to D to obtain $D_1$. On the other hand, suppose that only one of the cycles has a negative value, say $C_1$. Then we choose whichever of $C_1$ and $P_3$ has the smaller negative value.
We follow the same procedure starting off with each value i of MIN(M)(i, j) for j = 2, 3, …[log n]+1.
Once we have obtained the smallest (negative) cycle, say C*, from among all cycles tested, we multiply it by D to obtain DC* = $D_1$, a derangement, i.e., *a permutation which moves every point of* V = {1,2,…,n}. We continue with this procedure until we no longer can obtain a negative cycle. At this stage, we obtain *all negative* entries from MIN(M)(i, j) (i = 1,2, …,[log n] + 1; j= 1,2, … [log n] + 1}. We then apply the algorithm given above to each entry. If we obtain no further negative cycle, we go on to Phase 2.

## PHASE 2.

Phase 2 uses a modified form of the Floyd-Warshall Algorithm to obtain an optimal solution to the Assignment Problem. For short, we call this algorithm *the F-W n-vs algorithm,* i.e., the Floyd-Warshall negatively valued subpath algorithm. Again, we mention here that we use the version of the algorithm given in [6]. Assume that we have obtained $D_i^{-1}M^-$ using Phase 1. Given an entry, d(i, j), of a distance matrix, the basic idea of the F-W algorithm, is to see if there exists a path of length 2 (a triangular path) of form d(i, k) + d(k, j) the sum of whose values is less than d(i, j). If so, we replace d(i, j) by the value d(i, k) + d(k, j). For each value of j, we construct table of n columns and n rows in which we place k in the entry (i, j). We thus can keep track of both the new value of d(i, j) as well as the path of length two replacing it. (Note that as we go along the first "arc" generally represents a path.) We will use both MIN(M) as well as a balanced, binary search tree, T, containing all negative values of $D_i^{-1}M^-$. Each column heads a branch, while the rows numbers (each corresponding to the value of a negative path) are arranged in increasing order of magnitude. The value of the negative path at entry (a, b) is included at (a, b) on T. T is used in the following case: We've reached a point in which we can obtain no negative cycle, C, in $D_i^{-1}M^-(n)$ even though one exists. We then must continue using our algorithm, representing it now in the form $D_i^{-1}M^-(kn)(k = 2, 3, \dots)$. The maximum value that k could take is n-1 since each iteration going from k to k+1 adds one arc to any negative cycle obtained. (Note that an *iteration* is using the F-W n-vs algorithm through all columns from 1 through n in $D_i^{-1}M^-(in)$ when i > 1. ) As we go along, if a path no longer can be extended, we write its value on T in italics. The same is true once we've done all extensions of a path out of (a, b): We write the value of the path at (a, b) in italics. If $D_i$ contains only negatively-valued arcs, by construction, every arc is contained in a negatively-valued subpath of $D_i^{-1}M^-(n)$, the value matrix obtained after applying the F-W n-vs algorithm to every column of $D_i^{-1}M^-$. As an example of why i may be greater than 1, what may occur is that $C = \sum_{i=1}^{r} N_i + \sum_{j=1}^{s} P_i$ where each subpath $N_i$ is negatively-valued and always has an initial arc which has a negative value, while each arc of every $P_j$ always has a positive value. Since we cannot create such subpaths using the F-W n-vs algorithm, after we have completed using it in $D_i^{-1}M^-(n)$, the latter matrix may contain a set of disjoint negatively-valued subpaths which cannot be connected. The following theorem deals with that possibility. We first make a couple of definitions. Let M >0, $\varepsilon > 0$, be real positive numbers. Assume that the smallest value that any of the negatively-valued subpaths $N_i$ has is $- M - \varepsilon$. An F-W M-valued algorithm is one in which we may initiate an iteration using an entry whose value is no greater than M. Further assume that no positive arc has a smaller value than $\varepsilon$. Before going on, a *determining vertex* of a cycle of at most value L is an initial vertex, d, of a path P that traverses C such that every subpath is of value at most L. From the corollary to Theorem 1, at least one such vertex, d, always exists in a cycle C of value M.

THEOREM 5. Let P be a negatively-valued path in $D_h^{-1}M^-(k)$, say

**P = [a, ... , i, ... , d, ..., i]. Let i be the first vertex in P which changes P from a simple path to a non-simple one, i.e., it is the first vertex which is repeated. Let C be the subpath of P forming the cycle (d .... i). Here d is a determining vertex of the negatively-valued cycle, C. Then we can always obtain C as an independent negative cycle in at most the next iteration, $D_h^{-1}M^-(k+1)$., i.e., using at most n more columns in the algorithm.**

**PROOF. The simplest way to prove this theorem is to show that it holds in the worst possible case. Let r be the number of columns required to obtain all of the arcs in the subpath A = [a, ... , i], s, the number of columns passed through in order to obtain all of the arcs in I = [i, .... , d], t, the number of columns necessary to obtain all of the arcs in the subpath D = [d, e ... ,i]. Let P = $A \bigcup I \bigcup D$. Assume the following: r + s + d $\leq$ n, all of the arcs in D other than (d, e) have non-negative values, d is a determining vertex of D, d = n, while each subsequent vertex in D is smaller than the previous one. The following example illustrates a construction for P when n = 20.**

**Example. P = $[1^{-20}3^{5}7^{-5}13^{12}15^{1}19^{3}20^{-18}18^{1}14^{1}6^{3}7]$,**

**C = (20 18 14 6 7 13 15 19).**

**Note that there are no inversions in the natural order of the numbers until we reach d = 20. From that point the numbers decrease in value until they reach i = 7. From our algorithm, we can only choose an initial entry which is negative.**

**d(1, 3) = -20, d(1, 3, 7) = -15, d(1, 3, 7, 13) = -20,**
**d(1, 3, 7, 13, 15) = -8, d(1, 3, 7, 13, 15, 19) = -7,**
**d(1, 3, 7, 13, 15, 19, 20) = -4. In the triangular inequalities of the algorithm, the above paths would correspond to the respective distances obtained: d(1, 3), d(1, 7), d(1, 13), d(1, 15), d(1, 19), d(1, 20). Since we obtain a sequence of negative arcs with no inversions, we can form the negative path A + I by the time we reach the column j = 20. On the other hand, starting C with the pair of arcs (20, 18), (18, 14), we obtain (20, 14) where d(20, 14) = -17, we can proceed no further until we reach column j = 14 in $D_h^{-1}M^-(k+1)$ where we obtain (20, 14)(14, 6) = (20, 6) with d(20, 6) = -16. The subsequent arcs in C are (6, 7), (7, 13), (13, 15), (15, 19), (19, 20). Thus, we can obtain the negative cycle C at j = 20 in $D_i^{-1}M^-(k+1)$. On the other hand, P reaches j = 7 in $D_h^{-1}M^-(k+1)$. Thus, we obtained C as a negative cycle independent of P less than n = 20 columns after P reached i. A more abstract proof could be presented, but I think, in this case, the example given is more transparent. The crucial idea is that P must pass through the subpaths I and D but in a different order than the way C was obtained. It is just as likely that C can be obtained independently *before* P becomes non-simple if the roles of I and D are reversed, i.e., A and I have a large number of inversions while D has none.**

**PHASE 3.**

**STEP 1. We use the last matrix obtained - namely, the one which contains no negative cycles - at the start of the algorithm. Our goal is to obtain an optimal tour of M, $\sigma_{TSPOPT}$.**

**STEP 2. Starting with the derangements obtained during Phase 2, we check to see if any of them are tours. We also check to see if any of the permutations generated in obtaining $D_i$ (i=r,r-1,r-2,...,1) yields a tour when multiplied by $D_i$. Here $D_r = \sigma_{APOPT}$. If none does, we repeat Phase 1 n(log n) times and choose the smallest tour (if any) obtained. There always must be at least one smallest tour, $\tau_0$, since we start Phase 1 with a randomly chosen tour. In the worked-out example 1, we've used the difference of the total distances of $\sigma_{APOPT}$ and $\sigma_1$, the tour of smallest distance obtained using the methods just discussed. $|\sigma_1| = 161$. Thus, $m_0 = 6$. Using the corollary to Theorem 1, we can assume that each value in $D_3^{-1}M^-$ that belongs to a cycle is no larger than $m_0 - 1$.**

**STEP 3(a). Given a tour with an initial upper bound, $m_0$, for $|\sigma_{TSPOPT}|$, we use the F-W n-nvs (*Floyd-Warshall non-negatively-valued subpath*) algorithm to obtain all lowest-valued non-negative cycles obtainable. Since we no longer have to worry about negative cycles, non-negative cycles are obtainable provided $|\sigma_{\tau_0}| \neq |\sigma_{TSPOPT}|$. We continually test for tours using sets of disjoint cycles obtained, to see if we can find a permutation, s, such that $\tau_0 s = \tau_1$ where the value of $\tau_1$, is of value less than $m_0$, say $m_1$. We then assume that the upper bound for each arc is no greater than $m_1$ as we continue the algorithm. We thus discard cycles of value not less than $m_i$ (i = 1, 2, ,,,). We continue the algorithm through at most n-1 iterations.**

**THEOREM 6. Suppose that $m_j$ is the last upper bound obtained in Step 3(a). Suppose that $\sigma_{TSPOPT} = \sigma_{APOPT}s$, let C = (a $a_1$ $a_2$ ... $a_s$ b) be an arbitrary disjoint cycle of s where a is a determining vertex of C. Then there always exists a cycle of value less than $m_j$ obtained in Step 3(a) that is of the form C' =(a $b_1$ $b_2$ ... $b_r$ b).**

**COROLLARY. If we can obtain no cycle in Step 3(a), then $\sigma_{TSPOPT} = \sigma_{FWTSPOPT}$.**

**PROOF. $\sigma_{TSPOPT} = \sigma_{APOPT}s$. Here s is a permutation in $S_n$. The value of s is less than $m_j$ since were assuming that $|\sigma_{TSPOPT}| < |\sigma_{FWTSPOPT}|$. Let C = (a $a_1$ $a_2$ ... $a_s$ b) be a disjoint cycle of s. Let d(a, $a_1$, $a_2$,...,$a_s$,b) = m and m U d(b, a) = m' $\leq m_j$. Assume that a is a determining vertex of C. Then there must exist a shortest path, P, in $D_i^{-1}M^-$ from a to b of value no greater than m. Therefore, $|P|$ + d(b, a) $\leq$ m' implying that there must be a cycle, C', obtainable in Step 3(a). The corollary follows directly from the fact that if we could obtain a cycle [a ... b] in Step 3(b), then a corresponding cycle is obtainable in Step 3(a).**

**Comment. Theorem 6 narrows the search for cycles in Step 3(b). Our use of determining vertices is limited to vertices in cycles obtained in Step 3(a) that are initial vertices of arcs whose values are smaller than $m_j$.**

**Theorem 7. [2] The average number of cycles in a permutation in $S_n$ is log n.**

**Corollary, The average number of vertices in a cycle of a permutation in $S_n$ is $\frac{n}{\log n}$.**

**STEP 3(b).**

**Let m* be the lowest bound obtained in Step 3(a) for the value of a cycle. The best approach to an exhaustive search for *all* possible non-negative cycles is to construct a set of tree-like structures, denoted by c-trees, each of which has as its root an initial vertex of arc of value less than m*. Furthermore, we have a strict upper bound for the value of any subpath obtained along a branch. Let the root of one of the c-trees we are constructing be $v_i$. After each addition of an arc, (a, b), to a branch, we check to see if adding the arc (b, $v_i$) yields a cycle of value less than m*. If not, we continue the process. From Theorem 6, for large n, the number of cycles in an arbitrary permutation approaches log n. It follows that the average number of points in an arbitrary cycle of a permutation in $S_n$ is approximately $[\frac{n}{\log n}]$. Beginning with the $[\frac{n}{\log n}]$-th node, we backtrack along any path that hasn't been ended to see if it contains a repeated vertex. If it does, we end the path. A branch also ends when there are no arcs of value less than m* that can extend the branch. In general, we concentrate on the *last value* of each path [a, ..., b] and check if the value of [a, ... , b, a] is less than m*. In general, if a cycle exists, it will occur at an earlier node of the c-tree as a simple path rather than as a non-simple path containing a subpath that is a cycle. If we obtain a cycle, we test to see if it yields a smaller tour. If not we continue. We collect all cycles and each time a new one is obtained, we check all sets of disjoint cycles containing it to see if the product is a permutation which when multiplied by $\sigma_{APOPT}$ yields a tour of value less than $|\sigma_{FWTSPOPT}|$. It may require as many as n-1 iterations to obtain all cycles. If no product improves the value of $\sigma_{FWTSPOPT}$, then $\sigma_{TSPOPT} = \sigma_{FWTSPOPT}$.**

**Comment. From the corollary to Theorem 1, we know that we need only choose as an initial value, i, an entry whose value is no greater than the sum, i, of the values in the cycle. As we noted earlier, one of the advantages of our procedure for finding $\sigma_{APOPT}$ is that one or more of the derangements obtained may be an n-cycle. In particular, if we repeated the first part of Phase 1 n(log n) times, using a random n-cycle each time, the probability of obtaining at least one n-cycle is**

$$1 - \left(1 - \frac{e}{n}\right)^{n(\log n)} \approx 1 - \frac{1}{n^e} \to 1$$

**for large n. Of course, there is no guarantee that the n-cycle obtained would be close in value to $|\sigma_{TSPOPT}|$. The method given always obtains $\sigma_{TSPOPT}$ in (at most) Step 3(b) because there always exists *some* permutation, p, such that $(\sigma_{APOPT})p = \sigma_{TSPOPT}$. Whether this procedure can be done in every case in a finite amount of time is an open question.**

**Example 1.**

**PHASE 1.**

**Let D = (1 2 3 4 5 6 7 8), while M is the following distance matrix:**

**M**

| | 1 | 2 | 3 | 4 | 5 | 6 | 7 | 8 | |
|---|---|---|---|---|---|---|---|---|---|
| 1 | ∞ | 23 | 99 | 17 | 12 | 99 | 18 | 24 | 1 |
| 2 | 43 | ∞ | 2 | 73 | 15 | 100 | 53 | 28 | 2 |
| 3 | 1 | 84 | ∞ | 19 | 53 | 68 | 44 | 34 | 3 |
| 4 | 89 | 41 | 45 | ∞ | 40 | 71 | 79 | 51 | 4 |
| 5 | 83 | 62 | 94 | 88 | ∞ | 36 | 6 | 50 | 5 |
| 6 | 61 | 62 | 98 | 50 | 29 | ∞ | 52 | 40 | 6 |
| 7 | 50 | 21 | 53 | 68 | 39 | 26 | ∞ | 25 | 7 |
| 8 | 16 | 42 | 61 | 54 | 81 | 34 | 92 | ∞ | 8 |
| | 1 | 2 | 3 | 4 | 5 | 6 | 7 | 8 | |

**In the next table, MIN(M), we give the column numbers of the values in each row of M arranged in increasing order of magnitude.**

**MIN(M)**

| | 1 | 2 | 3 | 4 | 5 | 6 | 7 | 8 |
|---|---|---|---|---|---|---|---|---|
| 1 | 5 | 4 | 7 | 2 | 8 | 3 | 6 | 1 |
| 2 | 3 | 5 | 8 | 1 | 7 | 4 | 7 | 2 |
| 3 | 1 | 4 | 8 | 7 | 5 | 6 | 2 | 3 |
| 4 | 5 | 2 | 3 | 8 | 6 | 7 | 1 | 4 |
| 5 | 7 | 6 | 8 | 2 | 1 | 4 | 3 | 5 |
| 6 | 5 | 8 | 4 | 7 | 1 | 2 | 3 | 6 |
| 7 | 2 | 8 | 6 | 5 | 1 | 3 | 4 | 7 |
| 8 | 1 | 6 | 2 | 4 | 3 | 5 | 7 | 8 |

**MIN(M)(1, 1) = (1, 5): 12;  MIN(M)(2, 1) = (2, 3): 2;**
**MIN(M)(3, 1) = (3, 1):  1;  MIN(M)(4, 1) = (4, 5): 40;**
**MIN(M)(5, 1) = (5, 7):  6;  MIN(6, 1) = (6, 5): 29;**
**MIN(M)(7, 1) = (7, 2):  31;  MIN(8, 1) = (8, 1) = 16.**
**We now obtain the values of DIFF(i) (i = 1, 2, ..., 8).**
**DIFF(1) = d(1, 5)  -  d(1, 2) = -11.**
**DIFF(2) = d(2, 3)  -  d(2, 3) = 0.**
**DIFF(3) = d(3, 1)  -  d(3, 4) = -18.**
**DIFF(4) = d(4. 5)  -  d(4, 5) = 0.**
**DIFF(5) = d(5, 7)  -  d(5, 6) = -30.**

**DIFF(6) = d(6, 5) - d(6, 7) = -23.**
**DIFF(7) = d(7, 2) - d(7, 8) = -4.**
**DIFF(8) = d(8, 1) - d(8, 1) = 0.**
**min { DIFF(i) | i = 1, 2, ..., 8} = DIFF(5) = -30.**
**We will start Phase 1 with 5.**
**Before doing so, we write D and $D^{-1}$ in ROW FORM.**

| | | | | | | | | |
|---|---|---|---|---|---|---|---|---|
| | -11 | 0 | -18 | 0 | -30 | -23 | -4 | 0 |
| | 1 | 2 | 3 | 4 | 5 | 6 | 7 | 8 |
| D = | | | | | | | | |
| | 2 | 3 | 4 | 5 | 6 | 7 | 8 | 1 |
| | 1 | 2 | 3 | 4 | 5 | 6 | 7 | 8 |
| $D^{-1}$ = | | | | | | | | |
| | 8 | 1 | 2 | 3 | 4 | 5 | 6 | 7 |

**TRIAL 1.**

$(5, 7) \to (5, 6) \quad -30$

$(6, 5) \to (6, 4) \quad -23$

$(4, 5) \to (4, 4)$

**(4, 5) is an arc of D.**

$(4, \text{MIN(M)}(4, 2)) = (4, 2) \qquad d(4, 2) - d(4, 5) = 41 - 40 = 1$

$(4, 2) \to (4, 1)$

$(1, 5) \to (1, 4) \quad -11$

**P = [5, 6, 4, 1, 4].**
**(1, 5) is obtained from the arc (1, 6) in M.**
**d(1, 6) - d(1, 2) = 99 - 23 = 76.**
**(4, 5) is obtained from (4, 6).**
**d(4, 6) - d(4, 5) = 71 - 40 = 31**
**(6, 5) is obtained from the loop (6, 6).**
$s_{11} = (5^{-30}6^{-23}4^{1}1^{76})$, $s_{112} = (5^{-30}6^{-23}4^{31})$.
$s_{113}$ **= (5 6) is not admissible since** $Ds_{113}$ **is not a derangement.**
**The value of** $s_{11}$ **is -22.**
$s_{12}$ **= (5 6)(4 1). As shown earlier, (5 6) is not admissible.**
**(1, 4) is obtained from the arc (1, 5) of M.**
**d(1, 5) - d(1, 2) = 12 - 23 = -11.**
**Thus, since the second cycle has a negative value, we rename it** $s_{12} = (4^{1}1^{-11})$**. A reasonable question is why we go to the bother of saving it:** $s_{11}$ **has a smaller value than** $s_{12}$**. Our reason for doing so is that in Phase 3, we want to obtain an 8-cycle with as small as possible value. We may thus want to check to see if** $Ds_{11}$ **is an 8-cycle. In any event,** $s_{112}$ **= (5 6 4) has the smallest negative value so far: -22.**
**TRIAL 2. MIN(M)(5, 2) = (5, 6).**

**(5, 6) → (5, 5)**

**(5, 6) is an arc of D. We thus can go no farther with 5.**

**Therefore, $s_1$ = (5 6 4). $Ds_1$ = $D_1$. We obtained $D_1$ by exchanging those rows of D that have as initial arcs 5, 6, 4, with rows (5, 7), (6, 5), (4, 6), respectively. We must construct DIFF values for these arcs. All other DIFF values remain the same.**

| | | | | | | | | |
|---|---|---|---|---|---|---|---|---|
| | -11 | 0 | -18 | -31 | 0 | 0 | -4 | 0 |
| | 1 | 2 | 3 | 4 | 5 | 6 | 7 | 8 |
| $D_1$ = | | | | | | | | |
| | 2 | 3 | 4 | 6 | 7 | 5 | 8 | 1 |
| | 1 | 2 | 3 | 4 | 5 | 6 | 7 | 8 |
| $D_1^{-1}$ = | | | | | | | | |
| | 8 | 1 | 2 | 3 | 6 | 4 | 5 | 7 |

**d(5, 7) - d(5, 7) = 0, d(6, 5) - d(6, 5) = 0, d(4, 5) - d(4, 6) = -31.**

**We start with 4.**

**TRIAL 1. MIN(M)(4, 1) = (4, 5).**

**(4, 5) → (4, 6) -31**

**(6, 5) → (6, 6)**

**(6, 5) is an arc of $D_1$.**

**(6, MIN(6, 2)) = (6, 8) d(6, 8) - d(6, 5) = 11.**

**(6, 8) → (6, 7) 11**

**(7, 2) → (7, 1) -4**

**(1, 5) → (1, 6) -11**

**P = [4, 6, 7, 1, 6].**

**The arc (1, 4) is derived from the arc (1, 6) of M.**

**d(1, 6) - d(1, 2) = 99 - 23 = 76.**

**d(7, 6) - d(7, 8) = 26 - 25 = 1.**

**(6, 6) is a loop. Thus, the cycle (4 6) is not an admissible cycle since $D_1(1, 6)$ is not a derangement.**

**$s_{11}$ = $(4^{-31}6^{11}7^{-4}1^{76})$, $s_{111}$ = $(4^{-31}6^{11}7^{1})$.**

**The value of $s_{11}$ is non-negative. The value of $s_{111}$ is -19.**

**The only possibility for $s_{21}$ is $(6^{11}7^{-4}1^{-11})$. The value of the latter cycle is -4. We keep $s_{21}$ = (6 7 1) in our bag of negative cycles.**

**TRIAL 2. MIN(M)(4, 2) = (4, 2). d(4, 2) - d(4, 6) = 41 - 71 = -30.**

**(4, 2) → (4, 1) -30**

**(1, 5) → (1, 6) -11**

**(6, 5) → (6, 6)**

**(6, 5) is an arc of $D_1$.**

**(6, MIN(M)(6, 2)) = (6, 8) d(6, 8) - d(6, 5) = 40 - 29 = 11.**

**(6, 8) → (6, 7) 11**

**(7, 2) → (7, 1) -4**

**P = [4, 1, 6, 7, 1].**

**(7, 4) is derived from (7, 6) in M.**

**d(7, 6) - d(7, 8) = 26 - 25 = 1.**

**(6, 4) is derived from the loop (6, 6).**
**(1, 4) is derived from (1, 6).**
**d(1, 6) - d(1, 2) = 99 - 23 = 76.**
$s_{21} = (4^{-30}1^{-11}6^{11}7^{1})$, $s_{211} = (4^{-30}1^{76})$.
$s_{21}$ **has a negative value: -29.** $s_{211}$ **has a non-negative value.**
$s_{21}$ **has the smallest negative value thus far.**
**(7, 1) is derived from the arc (7, 2) of M.**
**d(7, 2) - d(7, 8) = 21 - 25 = -4.**
**Then** $s_{22} = (1^{-11}6^{11}7^{-4})$**. The value of** $s_{22}$ **is -4. We therefore keep it in our bag of negative cycles.**
**TRIAL 3. MIN(M)(4, 3) = (4, 3). d(4, 3) - d(4, 6) = 45 - 71 = -26.**

$(4, 3) \rightarrow (4, 2)$ -26
$(2, 3) \rightarrow (2, 2)$.

**(2, 3) is an arc of** $D_1$.
**(2, MIN(M)(2, 2)) = (2, 5). d(2, 5) - d(2, 3) = 15 - 2 = 13.**

$(2, 5) \rightarrow (2, 6)$ 13
$(6, 5) \rightarrow (6, 6)$.

**(6, 5) is an arc of** $D_1$.

**(6, MIN(M)(6, 2)) = (6, 8). d(6, 8) - d(6, 5) = 11**
$(6, 8) \rightarrow (6, 7)$ 11
$(7, 2) \rightarrow (7, 1)$ -4
$(1, 5) \rightarrow (1, 6)$

**P = [4, 2, 6, 7, 1, 6].**
**(1, 4) is derived from (1, 6). d(1, 6) - d(1, 2) = 76.**
**d(7, 6) - d(7, 8) = 1**
**(6, 6) is a loop.**
**d(2, 6) - d(2, 3) = 100 - 2 = 98.**
$P = [4^{-26}2^{13}6^{11}7^{-4}1^{76}]$
**Perusing the values we obtained above, the only possible negative cycle obtainable is** $s_{313}$ **= (4 2 6 7) which has a value of -1.**
$s_{32} = (6\ 7\ 1) = s_{22}$**; its value is -1.**
**Thus,** $s_2 = s_{21}$ **= (4 1 6 7).**
**In order to obtain** $D_2$**, we respectively exchange rows (4, 6), (1, 2), (6, 5), (7, 8) of** $D_1$ **with rows (4, 2), (1, 5), (6, 8), (7, 6).**
**Our new rows have the following DIFF values:**
**d(4, 5) - d(4, 2) = -1, d(1, 5) - d(1, 5) = 0, d(6, 5) - d(6, 8) = -11, d(7, 2) - d(7, 6) = -5.**

| | | | | | | | | |
|---|---|---|---|---|---|---|---|---|
| | 0 | 0 | -18 | -1 | 0 | -11 | -5 | 0 |
| | 1 | 2 | 3 | 4 | 5 | 6 | 7 | 8 |
| $D_2$ = | | | | | | | | |
| | 5 | 3 | 4 | 2 | 7 | 8 | 6 | 1 |

$$D_2^{-1} = \begin{matrix} 1 & 2 & 3 & 4 & 5 & 6 & 7 & 8 \\ 8 & 4 & 2 & 3 & 1 & 7 & 5 & 6 \end{matrix}$$

**We start with 3.**
**TRIAL 1.**
**(3, 1) $\rightarrow$ (3, 18) -18**
**(8, 1) $\rightarrow$ (8, 8)**
**(8, 1) is an arc of $D_2$.**
**(8, MIN(M)(8, 2)) = (8, 6). d(8, 6) - d(8, 1) = 18.**
**The path is non-negative. If we can go no further using the smallest negative DIFF, we are allowed to apply the procedure to at most the next [log n] smallest negative DIFF values.**
**Thus, we start with 6.**
**TRIAL 1.**
**(6, 5) $\rightarrow$ (6, 1) -11**
**(1, 5) $\rightarrow$ (1, 1)**
**(1, 5) is an arc of $D_2$.**
**(1, MIN(M)(1, 2)) = (1, 4). d(1, 4) - d(1,5) = 5.**
**(1, 4) $\rightarrow$ (1, 3) 5**
**(3, 1) $\rightarrow$ (3, 8) -18**
**(8, 1) $\rightarrow$ (8, 8)**
**(8, 1) is an arc of $D_2$.**
**(8, MIN(M)(8, 2)) = (8, 6). d(8, 6) - d(8, 1) = 18.**
**(8, 6) $\rightarrow$ (8, 7) 18**
**(7, 2) $\rightarrow$ (7, 4) -5**
**(4, 5) $\rightarrow$ (4, 1)**
**P = [6 1 3 8 7 4 1].**
**(4, 6) is derived from (4, 8).**
**d(4, 8) - d(4, 2) = 51 - 41 = 10.**
**d(7, 8) - d(7, 6) = 25 - 26 = -1.**
**(8, 8) is a loop.**
**d(3, 8) - d(3, 4) = 34 - 19 = 15.**
**d(1, 8) - d(1, 5) = 24 - 12 = 12.**
**We thus obtain the following cycles from P:**
$s_{11} = (6^{-11}1^{5}3^{-18}8^{18}7^{-5}4^{10})$, $s_{111} = (6^{-11}1^{5}3^{-18}8^{18}7^{-1})$, $s_{113} = (6^{-11}1^{5}3^{15})$, $s_{114} = (6^{-11}1^{12})$.
**$s_{11}$ has a value of -1. $s_{111}$ has a value of -7. The other two cycles have non-negative values.**
**(4, 1) is derived from the arc (4, 5) in M.**
**d(4, 5) - d(4, 2) = -1.**
**$s_{12} = (1^{5}3^{-18}8^{18}7^{-5}4^{-1})$. The value of $s_{12}$ is -1.**
**TRIAL 2. MIN(M)(6, 2)) = (6, 8). d(6, 8) - d(6, 8) = 0.**
**We can proceed no further using 6.**
**Since we have obtained a negative cycle using 6, we need not go to 7.**

**Thus, let $s_3 = s_{111} = (6\ \ 1\ \ 3\ \ 8\ \ 7)$. $D_3 = D_2 s_3$.**

**Rows (6, 8), (1, 5), (3, 4), *, 1), (7, 6) are replaced, respectively, by rows (6, 5), (1, 4), (3, 1), (8, 6), 7, 8) to obtain $D_3$.**

**The following are the new DIFF values:**
**d(6, 5) - d(6, 5) = 0.**
**d(1, 5) - d(1, 4) = -5.**
**d(3, 1) - d(3, 1) = 0.**
**d(8, 1) - d(8, 6) = -18.**
**d(7, 2) - d(7, 8) = -4.**

$$\begin{matrix} & -5 & 0 & 0 & -1 & 0 & 0 & -4 & -15 \\ & 1 & 2 & 3 & 4 & 5 & 6 & 7 & 8 \\ D_3 = & & & & & & & & \\ & 4 & 3 & 1 & 2 & 7 & 5 & 8 & 6 \end{matrix}$$

$$\begin{matrix} & 1 & 2 & 3 & 4 & 5 & 6 & 7 & 8 \\ D_3^{-1} = & & & & & & & & \\ & 3 & 4 & 2 & 1 & 6 & 8 & 5 & 7 \end{matrix}$$

**We now start with 8.**
**TRIAL 1.**
**(8, 1) $\to$ (8, 3) -18**
**(3, 1) $\to$ (3, 3).**
**(3, 1) is an arc of $D_3$.**
**(3, MIN(M)(3, 1) = (3, 4).**
**(3, 4) $\to$ (3, 1) 18**
**The path P = [8, 3, 1] is non-negative.**
**TRIAL 2. MIN(M)(8, 2) = (8, 6). d(8, 6) - d(8, 6) = 0. Thus, (8, 6) is an arc of $D_3$. We can go no further with 8.**

**We now start with 1.**
**TRIAL 1.**
**(1, 5) $\to$ (1, 6) -5**
**(6, 5) $\to$ (6, 6).**
**(6, 5) is an arc of $D_3$.**
**(6, MIN(M)(6, 2)) = (6, 8). d(6, 8) - d(6, 5) = 11.**
**P = [1, 6, 8] is non-negative.**
**(6, 1) is derived from (6, 4) of M. d(6, 4) - d(6, 5) = 50 - 29 = 21.**
**The cycle (1 6) has a non-negative value.**
**TRIAL 2. MIN(M)(1, 2) = (1, 4). d(1, 4) - d(1, 4) = 0.**
**(1, 4) is an arc of $D_3$. We can go no farther with 1.**

**The third smallest negative DIFF value is -4. It occurs at 7.**
**We start with 7.**
**TRIAL 1.**
**(7, 2) $\to$ (7, 4) -4**
**(4, 5) $\to$ (4, 6) -1**

(6, 5) $\rightarrow$ (6, 6).

(6, 5) is an arc of $D_3$.

(6, MIN(M)(6, 2)) = (6, 8). d(6, 8) - d(6, 5) = 11.

We can go no further since P = [7, 4, 6, 8] is non-negative.

(6, 7) is derived from (6, 8) of M.

d(6, 8) - d(6, 5) = 11.

d(4, 8) - d(4, 5) = 51 - 40 = 11.

$s_{11} = (7^{-4}4^{-1}6^{11})$, $s_{111} = (7^{-4}4^{11})$.

Both cycles have non-negative values.

TRIAL 2. MIN(M)(7, 2)) = (7, 8). d(7, 8) - d(7, 8) = 0.

(7, 8) is an arc of $D_3$. We thus conclude Phase 1 with

$D_3 = (1^{17}4^{41}2^{2}3^{1})(5^{6}7^{25}8^{34}6^{29})$. We next go to Phase 2 to try to obtain a derangement with a smaller value than $D_3$.

## PHASE 2.

$D_3^{-1}M$

| | 4 | 3 | 1 | 2 | 7 | 5 | 8 | 6 |
|---|---|---|---|---|---|---|---|---|
| | 1 | 2 | 3 | 4 | 5 | 6 | 7 | 8 |
| 1 | 17 | 99 | ∞ | 23 | 18 | 12 | 24 | 99 |
| 2 | 73 | 2 | 43 | ∞ | 53 | 15 | 28 | 100 |
| 3 | 19 | ∞ | 1 | 84 | 44 | 53 | 34 | 68 |
| 4 | ∞ | 45 | 89 | 41 | 79 | 40 | 51 | 71 |
| 5 | 88 | 94 | 83 | 62 | 6 | ∞ | 50 | 36 |
| 6 | 50 | 98 | 61 | 62 | 52 | 29 | 40 | ∞ |
| 7 | 68 | 53 | 50 | 21 | ∞ | 39 | 25 | 26 |
| 8 | 54 | 61 | 16 | 42 | 92 | 81 | ∞ | 34 |

$D_3^{-1}M^{-}$

| | 4 | 3 | 1 | 2 | 7 | 5 | 8 | 6 |
|---|---|---|---|---|---|---|---|---|
| | 1 | 2 | 3 | 4 | 5 | 6 | 7 | 8 |
| 1 | 0 | 82 | ∞ | 6 | 1 | - 5 | 7 | 82 |
| 2 | 71 | 0 | 41 | ∞ | 51 | 13 | 26 | 98 |
| 3 | 18 | ∞ | 0 | 83 | 43 | 52 | 33 | 67 |
| 4 | ∞ | 4 | 48 | 0 | 38 | - 1 | 10 | 30 |
| 5 | 82 | 88 | 77 | 56 | 0 | ∞ | 44 | 30 |
| 6 | 21 | 69 | 32 | 33 | 23 | 0 | 11 | ∞ |
| 7 | 43 | 28 | 25 | - 4 | ∞ | 14 | 0 | 1 |
| 8 | 20 | 27 | - 18 | 8 | 58 | 47 | ∞ | 0 |

j = 1, 2.

There are no negative entries in the first and second columns.

j = 3.
The only negative entry in the third column is in row 8.
d(8, 3) = - 18. MIN(M)(3, 1) = 1. $\sigma_4^{-1}(1) = 3$. (3, 3) is a loop.
MIN(M)(3, 2) = 4. $\sigma_4^{-1}(4) = 1$. d(3, 1) = 18. d(8, 3) + d(3, 1) $\geq$ 0.
MIN(M)(3, i), i = 3,4,5,6,7,8 are each greater than 18. Thus,
d(8, 3) + d(3, i) > 0, i = 3,4,5,6,7,8.
j = 4.
The only negative entry in column four occurs in row 7.
d(7, 4) = - 4. MIN(M)(4, 1) = 5. $\sigma_4^{-1}(5) = 6$. d(4, 6) = - 1.
d(7, 4) + d(4, 6) = - 5. < d(7, 6). Thus, we place - 5 in (7, 6).
MIN(M)(4, 2) = 2. $\sigma_4^{-1}(2) = 4$. (4, 4) is a loop.
MIN(M)(4, 3) = 3. $\sigma_4^{-1}(3) = 2$. d(4, 2) = 4. d(7, 4) + d(4, 6) $\geq$ 0. The
MIN(M)(4, i) for i = 4,5,6,7,8 yield values in $\sigma_4^{-1}M^{-}$ greater than 4. Thus,
d(7, 4) + d(4, i) > 0 for i = 4,5,6,7,8.

$P_4(4)$

| | 1 | 2 | 3 | 4 | 5 | 6 | 7 | 8 |
|---|---|---|---|---|---|---|---|---|
| 1 | | | | | | | | |
| 2 | | | | | | | | |
| 3 | | | | | | | | |
| 4 | | | | | | | | |
| 5 | | | | | | | | |
| 6 | | | | | | | | |
| 7 | | | | | | 1 | | |
| 8 | | | | | | | | |

j = 5.
There are no negative entries in column 5.
j = 6.
There are three negative entries in column 6: rows 1, 4 and 7.
d(1, 6) = - 5. MIN(M)(6, 1) = 5. $\sigma_4^{-1}(5) = 6$. (6, 6) is a loop.
MIN(M)(6, 2) = 8. $\sigma_4^{-1}(8) = 7$. d(6, 7) = 11. d(1, 6) + d(6, 7) > 0.
For i = 3,4,5,6,7,8, MIN(6, i) > MIN(M)(6, 2). Thus,
d(1, 6) + d(6, i) > 0.
d(4, 6) = - 1. MIN(M)(6, 2) = 8. $\sigma_4^{-1}(8) = 7$. d(6, 7) = 11.
d(4, 6) + d(6, 7) > 0. The same is true with
MIN(M)(6, i), i = 3,4,5,6,7,8.
d(7, 6) = - 5. MIN(M)(6, 2) = 8. $\sigma_4^{-1}(8) = 7$. d(6, 7) = 11.
d(4, 6) + d(6, 7) > 0. The same is true with
MIN(M)(6, i), i = 3,4,5,6,7,8.
j = 7, 8.
There are no negative entries in columns 7 and 8.

We have only one negative path in $D_3^{-1}M^-(8)$: [7, 4, 6] with a value of -5. This path can't be extended in $D_3^{-1}M^-(16)$. It follows that $\sigma_{APOPT} = D_3$. The total sum of its values is 155.

## PHASE 3.

We now apply Step 3(a) to obtain an approximation to $\sigma_{TSPOPT}$.
In order to obtain the smallest possible initial tour, we go back to $D_2$ to check to see if any of the permutations we obtained yields an 8-cycle when multiplied by $D_2$, i.e., a tour. We first check those cycles that have negative values: $s_{11}$ and $s_{21}$ both have values of -1.
$D_2s_{11} = (1^{17}4^{51}8^{34}6^{29}5^{6}7^{21}2^{2}3^{1})$. Its value is 161 - 6 greater than $\sigma_{APOPT}$.
Since 6 is so small, we will let $m_0 = 6$. Furthermore, define $\sigma_1 = D_2s_{11}$.
We thus must be able to obtain a positive cycle of value no greater than 5 using only entries of $D_3^{-1}M^-$ no greater than 5. We use the F-W non-n-vs algorithm.
Comment. If we hadn't been able to obtain a tour using permutations obtained while working to obtain $D_3$, we would use those obtained trying to obtain $D_2$, etc. . As we go along, we may obtain cycles that lower the upper bound from $m_0 = 6$ to $m_1$. From that moment onward, we only consider paths whose value is less than $m_1$. This may occur a number of times during one iteration (j = 1 through j = 8). Therefore, the best approach is to complete the iteration and underline in $D_3^{-1}M^-(8)$ and $P_8$ only those entries whose values are less than the last $m_i$ obtained. Henceforth, we will only give each sum of two arcs and its respective value in obtaining new entries. We underline all paths that we are currently trying to extend. We place in italics all entries that we've passed through always using the *smallest* value obtained in a path.

$D_3^{-1}M^-$

| | 4 | 3 | 1 | 2 | 7 | 5 | 8 | 6 |
|---|---|---|---|---|---|---|---|---|
| | 1 | 2 | 3 | 4 | 5 | 6 | 7 | 8 |
| 1 | 0 | 82 | ∞ | 6 | 1 | - 5 | 7 | 82 |
| 2 | 71 | 0 | 41 | ∞ | 51 | 13 | 26 | 98 |
| 3 | 18 | ∞ | 0 | 83 | 43 | 52 | 33 | 67 |
| 4 | ∞ | 4 | 48 | 0 | 38 | - 1 | 10 | 30 |
| 5 | 82 | 88 | 77 | 56 | 0 | ∞ | 44 | 30 |
| 6 | 21 | 69 | 32 | 33 | 23 | 0 | 11 | ∞ |
| 7 | 43 | 28 | 25 | - 4 | ∞ | 14 | 0 | 1 |
| 8 | 20 | 27 | - 18 | 8 | 58 | 47 | ∞ | 0 |

j=3.
(8, 3)(3, 1): 0.
j=4.

**(7, 4)(4, 2): 0; (7, 4)(4, 6): -5.**
**j=8.**
**(7, 8)(8, 3): -17.**

**$D_3^{-1}M^-(8)$**

| | 4<br>1 | 3<br>2 | 1<br>3 | 2<br>4 | 7<br>5 | 5<br>6 | 8<br>7 | 6<br>8 |
|---|---|---|---|---|---|---|---|---|
| 1 | 0 | 82 | ∞ | 6 | 1 | - 5 | 7 | 82 |
| 2 | 71 | 0 | 41 | ∞ | 51 | 13 | 26 | 98 |
| 3 | 18 | ∞ | 0 | 83 | 43 | 52 | 33 | 67 |
| 4 | ∞ | 4 | 48 | 0 | 38 | - 1 | 10 | 30 |
| 5 | 82 | 88 | 77 | 56 | 0 | ∞ | 44 | 30 |
| 6 | 21 | 69 | 32 | 33 | 23 | 0 | 11 | ∞ |
| 7 | 43 | 0 | -17 | -4 | ∞ | -5 | 0 | 1 |
| 8 | 0 | 27 | - 18 | 8 | 58 | 47 | ∞ | 0 |

**j=1.**
**(8, 1)(1, 6): -5; (8, 1)(1, 5): 1.**
**j=3.**
**(7, 3)(3, 1): 1.**

**$D_3^{-1}M^-(16)$**

| | 4<br>1 | 3<br>2 | 1<br>3 | 2<br>4 | 7<br>5 | 5<br>6 | 8<br>7 | 6<br>8 |
|---|---|---|---|---|---|---|---|---|
| 1 | 0 | 82 | ∞ | 6 | 1 | - 5 | 7 | 82 |
| 2 | 71 | 0 | 41 | ∞ | 51 | 13 | 26 | 98 |
| 3 | 18 | ∞ | 0 | 83 | 43 | 52 | 33 | 67 |
| 4 | ∞ | 4 | 48 | 0 | 38 | - 1 | 10 | 30 |
| 5 | 82 | 88 | 77 | 56 | 0 | ∞ | 44 | 30 |
| 6 | 21 | 69 | 32 | 33 | 23 | 0 | 11 | ∞ |
| 7 | 1 | *0* | *-17* | -4 | ∞ | *-5* | 0 | 1 |
| 8 | *0* | 27 | - 18 | 8 | 1 | -5 | ∞ | 0 |

**j=1.**
**(7, 1)(1, 5): 2**
**(7, 5) cannot be extended to a path of value less than 6.**
**It follows that $\sigma_{FWTSPOPT} = \sigma_1 = (1^{17}4^{51}8^{34}6^{29}7^{21}2^{2}3^{1})$.**
**Furthermore, from the corollary to Theorem 6, since we are unable to obtain a positive cycle of value less than 6, $\sigma_{TSPOPT} = \sigma_1$.**